\documentclass[11pt]{article}

\usepackage{latexsym}
\usepackage{amssymb}

\newtheorem{theorem}{Theorem}

\newtheorem{corollary}{Corollary}

\newtheorem{lemma}{Lemma}

\begin{document}
{
\begin{center}
{\Large\bf
On the multiplication operator by an independent variable in matrix Sobolev spaces.}
\end{center}
\begin{center}
{\bf S.M. Zagorodnyuk}
\end{center}


V. N. Karazin Kharkiv National University \newline\indent
School of Mathematics and Computer Sciences \newline\indent
Department of Higher Mathematics and Informatics \newline\indent
Svobody Square 4, 61022, Kharkiv, Ukraine

Sergey.M.Zagorodnyuk@gmail.com; Sergey.M.Zagorodnyuk@univer.kharkov.ua

\vspace{1cm}

\noindent
\textbf{Abstract}. We study the operator $\mathcal{A}$ of multiplication by an independent variable in a matrix Sobolev space $W^2(M)$.
In the cases of finite measures on $[a,b]$ with $(2\times 2)$ and $(3\times 3)$ real continuous matrix weights of full rank it is shown that     
the operator $\mathcal{A}$ is symmetrizable. Namely, there exist two symmetric operators $\mathcal{B}$ and
$\mathcal{C}$ in a larger space such that $\mathcal{A} f = \mathcal{C} \mathcal{B}^{-1} f$,
$f\in D(\mathcal{A})$.
As a corollary, we obtain some new orthogonality conditions for the associated Sobolev orthogonal polynomials.
These conditions involve two symmetric operators in an indefinite metric space.

\noindent
\textbf{Keywords}: a multiplication operator, Sobolev space, a symmetrizable operator, Sobolev orthogonal polynomials. 

\noindent
\textbf{MSC 2010}: 47B37

\section{Introduction.}

The theory of Sobolev orthogonal polynomials has got a powerful impulse for development during the past 30 years, see the 
survey~\cite{cit_5150_M_X_Survey_2015} and references therein. 
However, many aspects of the theory are still hidden.
One of intriguing topics is an investigation of the associated multiplication operator by an independent variable,
which acts in the underlying Sobolev type space. Let us recall basic definitions.

Fix an arbitrary Borel subset $K$ of the complex plane and an arbitrary $\rho\in\mathbb{N}$. 
Let $M(\delta) = ( m_{k,l}(\delta) )_{k,l=0}^\rho$ be a $\mathbb{C}_{(\rho+1)\times (\rho+1)}^\geq$-valued function on 
$\mathfrak{B}(K)$, which entries are countably additive on $\mathfrak{B}(K)$ ($\delta\in\mathfrak{B}(K)$) (it is called \textit{a non-negative 
Hermitian-valued measure on $(K,\mathfrak{B}(K))$}).
Denote by $\tau(\delta)$ \textit{the trace measure}, $\tau(\delta) := \sum_{k=0}^\rho m_{k,k}(\delta)$, $\delta\in\mathfrak{B}(K)$.
By $M'_\tau := dM/d\tau = ( dm_{k,l}/ d\tau )_{k,l=0}^{\rho}$, we denote \textit{the trace derivative of $M$}, see~\cite{cit_7000__Rosenberg_1964}.
One means by $L^2(M)$ a set of all (classes of the equivalence of)
measurable vector-valued functions
$\vec f(z): K\rightarrow \mathbb{C}_{\rho+1}$, 
$\vec f = (f_0(z),f_1(z),\ldots,f_\rho(z))$, such that
$$ \| \vec f \|^2_{L^2(M)} := \int_K  \vec f(z) M'_\tau(z) \vec f^*(z) d\tau  < \infty. $$
It is known that $L^2(M)$ is a Hilbert space with the following scalar product:
\begin{equation}
\label{f2_3}
( \vec f, \vec g )_{L^2(M)} := \int_K  \vec f(z) M'_\tau(z) (\vec g(z))^* d\tau,\qquad  \vec f,\vec g\in L^2(M). 
\end{equation}
It is also known (see~\cite[p. 294]{cit_7000__Rosenberg_1964}, \cite[Lemma 2.1]{cit_8000__Weron_1974})
that one can consider an arbitrary $\sigma$-finite (non-negative) measure $\mu$, with respect to which all $m_{k,l}$ are
absolutely continuous, and set $M_0(z) = M_{0,\mu}(z) := dM/d\mu$ (the Radon-Nikodym derivative of $M$ with respect to $\mu$).
The integral in~(\ref{f2_3}) exists if and only if the following integral exists:
\begin{equation}
\label{f2_4}
\int_K  \vec f(z) M_{0,\mu}(z) (\vec g(z))^* d\mu. 
\end{equation}
If the integrals exist, they are equal. Such measures $\mu$ we shall call \textit{admissible (for $M$)}. The matrix function $M_0(z) = M_{0,\mu}(z)$
is said to be \textit{the (matrix) weight, corresponding to an admissible measure $\mu$}.
If $\det M_0(z) > 0$, $\forall z\in K$, we shall say that \textit{the matrix weight $M_0(z)$ has full rank}.

Denote by $A^2(M)$ a linear manifold in $L^2(M)$ including those classes of the equivalence $[\cdot]$ which possess 
a representative of the following form: 
\begin{equation}
\label{f2_5}
\vec f(z) = (f(z),f'(z),...,f^{(\rho)}(z)). 
\end{equation}
By $W^2(M)$ we denote the closure of $A^2(M)$ in the norm of $L^2(M)$. 
The subspace $W^2(M)$ is said to be \textit{the Sobolev space with the matrix measure $M$}, see~\cite{cit_10000__Z_JDEA_2021}.
Elements of $A^2(M)$ will be also
denoted by their first components.

Let $\mu$ be an admissible measure and $M_0(z)$ be the corresponding matrix weight.
Suppose that $1,z,z^2,...$, all belong to $W^2(M)$. 
Assume that
\begin{equation}
\label{f2_7}
(p,p)_{W^2(M)} > 0,
\end{equation}
for an arbitrary non-zero $p\in\mathbb{P}$.
Then one can apply the Gram-Schmidt orthogonalization process to construct a system
$\{ y_n(z) \}_{n=0}^\infty$, $\deg y_n=n$, of the associated Sobolev orthogonal polynomials:
$$ \int_K  (y_n(z),y_n'(z),...,y_n^{(\rho)}(z)) M_0(z) 
\overline{
\left(
\begin{array}{cccc} y_m(z)\\
y_m'(z)\\
\vdots\\
y_m^{(\rho)}(z)
\end{array}
\right)
}
d\mu
= 
$$
\begin{equation}
\label{f2_10}
= A_n \delta_{n,m},\qquad               A_n>0,\quad           n,m\in\mathbb{Z}_+.
\end{equation}

In this paper we shall only consider the case of $K = [a,b]$, $-\infty < a < b < +\infty$.
We fix an admissible measure $\mu$ and assume that the entries of the corresponding matrix weight $M_0(z)$ are real-valued continuous functions.
In this case the matrix weight $M_0(z)$ is said to be \textit{real-valued continuous}.
Denote by $C(K)$ the set of all classes of equivalence of functions in $L^2(M)$, which include a continuous representative on $[a,b]$
(i.e. each entry of the representative is continuous on $[a,b]$).
We denote by $\mathcal{A}$ the operator on the whole $A^2(M)$, which sends $[f(z)]$ to $[z f(z)]$.

The main purpose of this paper is to show that 
under some general conditions
there exist two symmetric operators $\mathcal{B}$ and
$\mathcal{C}$ in $L^2(M)$ ($D(B)=D(C)=C(K)$) such that $\mathcal{A}f = 
\mathcal{C} \mathcal{B}^{-1} f$, $f\in D(\mathcal{A})$.
As a corollary, we obtain some new orthogonality conditions for the associated Sobolev orthogonal polynomials.
These conditions involve two symmetric operators in an indefinite metric space.
Observe that some other interesting questions concerning matrix measures and Sobolev type inner products were studied
in~\cite{cit_2000__Kwon___2009}, \cite{cit_1500___Kim__2014}.

\noindent
{\bf Notations. }
As usual, we denote by $\mathbb{R}, \mathbb{C}, \mathbb{N}, \mathbb{Z}, \mathbb{Z}_+$,
the sets of real numbers, complex numbers, positive integers, integers and non-negative integers,
respectively. 
By $\mathbb{Z}_{k,l}$ we mean all integers $j$ satisfying the following inequality:
$k\leq j\leq l$; ($k,l\in\mathbb{Z}$).
By $\mathbb{C}_{m\times n}$ we denote the set of all $(m\times n)$ matrices with complex entries,
$\mathbb{C}_{n} := \mathbb{C}_{1\times n}$, 
$\mathbb{C}^{n} := \mathbb{C}_{n\times 1}$, $m,n\in\mathbb{N}$.
By $\mathbb{C}_{n\times n}^\geq$ we denote the set of all nonnegative Hermitian matrices from $\mathbb{C}_{n\times n}$,
$n\in\mathbb{N}$.
For $A\in\mathbb{C}_{m\times n}$ the notation $A^*$ stands for the adjoint matrix ($m,n\in\mathbb{N}$),
and $A^T$ means the transpose of $A$.
By $\mathbb{P}$ we denote the set of all polynomials with complex coefficients.
For an arbitrary Borel subset $K$ of the complex plane we denote by $\mathfrak{B}(K)$ the set of all Borel subsets of $K$.
Let $\mu$ be an arbitrary (non-negative) measure on $\mathfrak{B}(K)$.
By $L^2_\mu = L^2_{\mu,K}$ we denote the usual space of (the classes of the equivalence of) complex Borel measurable functions $f$ on $K$ such that
$\| f \|_{L^2_{\mu,K}}^2 := \int_K |f(z)|^2 d\mu < \infty$.

By
$(\cdot,\cdot)_H$ and $\| \cdot \|_H$ we denote the scalar product and the norm in a Hilbert space $H$,
respectively. The indices may be omitted in obvious cases.
For a set $M$ in $H$, by $\overline{M}$ we mean the closure of $M$ in the norm $\| \cdot \|_H$.

\section{Matrix multiplication operators and symmetrizability.}

Our main objective here is to deduce the following theorem.

\begin{theorem}
\label{t3_1}
Let $\rho\in\{ 1,2 \}$, $K=[a,b]$, $-\infty < a < b < +\infty$, and $M(\delta) = ( m_{k,l}(\delta) )_{k,l=0}^\rho$ be 
a non-negative Hermitian-valued measure on $(K,\mathfrak{B}(K))$. Suppose that $\mu$ is a finite admissible measure for $M$, such that
the corresponding matrix weight $M_0(z)$ is real-valued continuous and of full rank.
Then the operator $\mathcal{A}$:
\begin{equation}
\label{f3_5}
\mathcal{A} [f(z)] = [ zf(z) ],\qquad [f(z)]\in A^2(M),
\end{equation}
is well-defined, and it has the following representation:
\begin{equation}
\label{f3_7}
\mathcal{A} = \mathcal{C} \mathcal{B}^{-1} |_{A^2(M)}.
\end{equation}
Here $\mathcal{B},\mathcal{C}$ are some linear symmetric operators in $L^2(M)$, and $\mathcal{B}$ is invertible.

\end{theorem}

\noindent
\textbf{Proof.}
At first we shall consider an arbitrary $\rho\in\mathbb{N}$, and $K$, $M(\delta)$, $\mu$, $M_0(z)$ like in the statement of the theorem.
In order to construct the required operators $\mathcal{B},\mathcal{C}$, we shall use matrix multiplication operators in $L^2(M)$.

\begin{lemma}
\label{l3_1}
Let $\rho\in\mathbb{N}$, and $K$, $M(\delta)$, $\mu$, $M_0(z)$ be like in the statement of Theorem~\ref{t3_1}.
Let $D(z) = (d_{i,j}(z))_{i,j=0}^\rho$, be an arbitrary $(\rho+1)\times(\rho+1)$ matrix function, which entries are
real-valued continuous functions on $K$.
Then the following operator $\mathcal{D}$:
\begin{equation}
\label{f3_16}
\mathcal{D} u = [(f_0(z),...,f_\rho(z)) D(z)],\qquad u = [(f_0(z),...,f_\rho(z))]\in C(K),
\end{equation}
is a well-defined linear operator in $L^2(M)$ with the domain $C(K)$.
\end{lemma}

\noindent
\textbf{Proof of lemma.}
Suppose that an element $u\in C(K)$ has two continuous representatives $\vec f$, $\vec g$:
$$ \| \vec f - \vec g \|_{L^2(M)} = 0. $$
Observe that $(M_0(z))^{\frac{1}{2}}$ has $\mu$-measurable entries (cf.~\cite{cit_7000__Rosenberg_1964} for the case of the trace derivative).
In fact, the norm of the operator $M_0(z)$ in $\mathbb{C}^{\rho+1}$ is a continuous function  in $z$, and therefore it attains its maximum $L$ on $K$.
Thus, the spectra of all of $M_0(z)$ lie in $[0,L]$.
On the segment $[0,L]$ we can approximate $\sqrt{x}$ by a polynomial $p_k(x)$ in the uniform norm $\| \cdot \|_U$:
$$ \| \sqrt{x} - p_k(x) \|_U < \frac{1}{k},\qquad k\in\mathbb{N}. $$
Then
$$ \| (M_0(z))^{\frac{1}{2}} \vec x - p_k(M_0(z)) \vec x \|^2_{\mathbb{C}^{\rho+1}} = 
\int_0^L |\sqrt{\lambda} - p_k(\lambda)|^2 d(E_\lambda(z) \vec x,\vec x)\leq $$
$$ \leq \frac{1}{k^2} \| \vec x \|^2, $$
where $E_\lambda(z)$ is the orthogonal resolution of the identity for $B_0(z)$.
Thus, elements of $(M_0(z))^{\frac{1}{2}}$ are $\mu$-measurable, as they are the limits of $\mu$-measurable functions.

\noindent
By the structure of the inner product it follows that (cf.~\cite{cit_7000__Rosenberg_1964}):
$$ (\vec f - \vec g) M_0^{\frac{1}{2}}(z) = 0, $$
$\mu$-a.e. on $K$. Then
$$ (\vec f - \vec g) D(z)  =  (\vec f - \vec g) M_0(z) M_0^{-1}(z) D(z) =  0, $$
$\mu$-a.e. on $K$. Therefore
$$ \| \vec f D(z) - \vec g D(z) \|_{L^2(M)} = 0. $$
Consequently, the operator $\mathcal{D}$ is well-defined.
The linearity is obvious. The lemma is proved.

Let us return to the proof of the theorem.
By the Leibniz rule we may write:
$$ (z f(z))^{(r)} = z f^{(r)}(z) + r f^{(r-1)}(z),\qquad [f]\in A^2(M),\ 0\leq r\leq \rho. $$
Therefore
$$ (zf(z),(zf(z))',...,(zf(z))^{(\rho)}) =
(f(z),f'(z),...,f^{(\rho)}(z)) A(z), $$
where 
\begin{equation}
\label{f3_15}
A(z) := 
\left(
\begin{array}{cccccc}
z & 1 & 0 & \cdots & 0 & 0 \\
0 & z & 2 & \cdots & 0 & 0 \\
0 & 0 & z & \cdots & 0 & 0 \\
\vdots & \vdots & \vdots & \ddots & \vdots & \vdots \\
0 & 0 & 0 & \cdots & z & \rho \\
0 & 0 & 0 & \cdots & 0 & z \end{array}
\right).
\end{equation}

Denote by $\widehat{\mathcal{A}}$ the operator in $L^2(M)$, with $D(\widehat{\mathcal{A}}) = C(K)$, which maps
$[(f_0(z),...,f_\rho(z))]$ to $[(f_0(z),...,f_\rho(z)) A(z)]$.
Observe that $\widehat{\mathcal{A}} \supseteq \mathcal{A}$.
Thus the operator $\mathcal{A}$ is well-defined.
Let $B(z) = (b_{i,j}(z))_{i,j=0}^\rho$, be an arbitrary $(\rho+1)\times(\rho+1)$ matrix function, which entries are
real-valued continuous functions on $K$.
Set
\begin{equation}
\label{f3_17}
C(z) := B(z) A(z),\qquad z\in K. 
\end{equation}
By $\mathcal{B}$, $\mathcal{C}$ we mean the operators in $L^2(M)$, with $D(\mathcal{B}) = D(\mathcal{C}) = C(K)$, such that:
$$ \mathcal{B} [(f_0(z),...,f_\rho(z))] = [(f_0(z),...,f_\rho(z)) B(z)], $$
$$\mathcal{C} [(f_0(z),...,f_\rho(z))] = [(f_0(z),...,f_\rho(z)) C(z)]. $$
By~(\ref{f3_17}) we see that (notice the reversed order due to the right multiplication)
\begin{equation}
\label{f3_19}
\mathcal{C} = \mathcal{A} \mathcal{B}. 
\end{equation}

Our aim now is to specify functions $b_{i,j}(z)$ in order to get relation~(\ref{f3_7}), as well as to guarantee the symmetry of
$\mathcal{B}$ and $\mathcal{C}$.
The structure of the inner product in $L^2(M)$ shows that the following conditions provide the symmetry of $\mathcal{B}$ and $\mathcal{C}$:
\begin{equation}
\label{f3_20}
B(z) M_0(z) = M_0(z) B^*(z),\quad C(z) M_0(z) = M_0(z) C^*(z),\qquad z\in K. 
\end{equation}
Observe that
$$ C(z) = z B(z) + \widehat B(z),\quad \widehat B(z) := (j b_{i,j-1}(z))_{i,j=0}^\rho,\ z\in K, $$
where $b_{k,l}$ with negative indices are zeros. 
Thus conditions~(\ref{f3_20}) are satisfied, if and only if
\begin{equation}
\label{f3_22}
B(z) M_0(z) = M_0(z) B^*(z),\quad \widehat B(z) M_0(z) = M_0(z) \widehat B^*(z),\qquad z\in K. 
\end{equation}
Conditions~(\ref{f3_22}) are equivalent to the following conditions:
\begin{equation}
\label{f3_24}
B(z) M_0(z), \widehat B(z) M_0(z)\quad \mbox{are Hermitian matrices for all $z\in K$.} 
\end{equation}
Let $M_0(z) = (\widetilde m_{i,j})_{i,j=0}^\rho$.
Conditions~(\ref{f3_24}) in turn are equivalent to the following conditions:
\begin{equation}
\label{f3_26}
\sum_{k=0}^\rho b_{i,k} \widetilde m_{k,j}(z) = \sum_{k=0}^\rho b_{j,k} \widetilde m_{k,i}(z),\ 0\leq j\leq\rho-1,\ j<i\leq\rho;
\end{equation}
\begin{equation}
\label{f3_28}
\sum_{k=1}^\rho k b_{i,k-1} \widetilde m_{k,j}(z) = \sum_{k=1}^\rho k b_{j,k-1} \widetilde m_{k,i}(z),\ 0\leq j\leq\rho-1,\ j<i\leq\rho;\
z\in K.
\end{equation}

If one can find a real continuous solution $B(z)$ to equations~(\ref{f3_26}),(\ref{f3_28}), such that 
$$ \det B(z) \not=0,\qquad z\in K, $$
then the corresponding operators $\mathcal{B},\mathcal{C}$ will provide the required representation for $\mathcal{A}$.

\noindent
\textit{Case $\rho=1$}. Equations~(\ref{f3_26}),(\ref{f3_28}) now have the following form:
\begin{equation}
\label{f3_30}
\left\{
\begin{array}{cc}
b_{1,0} \widetilde m_{0,0}(z) + b_{1,1} \widetilde m_{1,0}(z) = b_{0,0} \widetilde m_{0,1}(z) + b_{0,1} \widetilde m_{1,1}(z),\\
b_{1,0} \widetilde m_{1,0}(z) = b_{0,0} \widetilde m_{1,1}(z),\end{array}
\right.\qquad z\in K.
\end{equation}
Then $b_{1,0}, b_{1,1}$ are arbitrary real continuous functions on $K$, while
$$ b_{0,0} = b_{1,0} \frac{ \widetilde m_{1,0} }{  \widetilde m_{1,1}  }, $$
$$ b_{0,1} = \frac{1}{ \widetilde m_{1,1} } \left(
b_{1,0} \left( 
\widetilde m_{0,0} - \frac{ \widetilde m_{1,0}^2 }{  \widetilde m_{1,1}  }
\right) + 
b_{1,1} \widetilde m_{1,0}
\right). $$
In particular, choose $b_{1,0} = \widetilde m_{1,1}, b_{1,1} = 0$, to get
$$ B(z) =
\left(
\begin{array}{cc}
\widetilde m_{1,0}(z) & \widetilde m_{0,0}(z) - \frac{ \widetilde m_{1,0}^2(z) }{  \widetilde m_{1,1}(z)  } \\
\widetilde m_{1,1}(z) & 0\end{array}
\right),\quad z\in K. $$
It is clear that $\det B(z) \not= 0$, $z\in K$.

\noindent
\textit{Case $\rho=2$}. Equations~(\ref{f3_26}),(\ref{f3_28}) lead to the following six equations:
\begin{equation}
\label{f3_44__1}
b_{1,0} \widetilde m_{0,0} + b_{1,1} \widetilde m_{1,0} + b_{1,2} \widetilde m_{2,0} = 
b_{0,0} \widetilde m_{0,1} + b_{0,1} \widetilde m_{1,1} + b_{0,2} \widetilde m_{2,1},
\end{equation}
\begin{equation}
\label{f3_44__2}
b_{1,0} \widetilde m_{1,0} + 2b_{1,1} \widetilde m_{2,0} = 
b_{0,0} \widetilde m_{1,1} + 2b_{0,1} \widetilde m_{2,1},
\end{equation}
\begin{equation}
\label{f3_44__3}
b_{2,0} \widetilde m_{0,0} + b_{2,1} \widetilde m_{1,0} + b_{2,2} \widetilde m_{2,0} = 
b_{0,0} \widetilde m_{0,2} + b_{0,1} \widetilde m_{1,2} + b_{0,2} \widetilde m_{2,2},
\end{equation}
\begin{equation}
\label{f3_44__4}
b_{2,0} \widetilde m_{1,0} + 2b_{2,1} \widetilde m_{2,0} = 
b_{0,0} \widetilde m_{1,2} + 2b_{0,1} \widetilde m_{2,2},
\end{equation}
\begin{equation}
\label{f3_44__5}
b_{2,0} \widetilde m_{0,1} + b_{2,1} \widetilde m_{1,1} + b_{2,2} \widetilde m_{2,1} = 
b_{1,0} \widetilde m_{0,2} + b_{1,1} \widetilde m_{1,2} + b_{1,2} \widetilde m_{2,2},
\end{equation}
\begin{equation}
\label{f3_44__6}
b_{2,0} \widetilde m_{1,1} + 2b_{2,1} \widetilde m_{2,1} = 
b_{1,0} \widetilde m_{1,2} + 2b_{1,1} \widetilde m_{2,2}.
\end{equation}
Equations~(\ref{f3_44__3})-(\ref{f3_44__6}) are equivalent to the following equations:
\begin{equation}
\label{f3_44__3_circle}
b_{0,2} = b_{2,0} \frac{\widetilde m_{0,0}}{\widetilde m_{2,2}} + b_{2,1} \frac{\widetilde m_{1,0}}{\widetilde m_{2,2}} + 
b_{2,2} \frac{\widetilde m_{2,0}}{\widetilde m_{2,2}} - 
b_{0,0} \frac{\widetilde m_{0,2}}{\widetilde m_{2,2}} - b_{0,1} \frac{\widetilde m_{1,2}}{\widetilde m_{2,2}},
\end{equation}
\begin{equation}
\label{f3_44__4_circle}
b_{0,1} = b_{2,0} \frac{\widetilde m_{1,0}}{2\widetilde m_{2,2}} + b_{2,1} \frac{\widetilde m_{2,0}}{\widetilde m_{2,2}} - 
b_{0,0} \frac{\widetilde m_{1,2}}{2\widetilde m_{2,2}},
\end{equation}
\begin{equation}
\label{f3_44__5_circle}
b_{1,2} = b_{2,0} \frac{\widetilde m_{0,1}}{\widetilde m_{2,2}} + b_{2,1} \frac{\widetilde m_{1,1}}{\widetilde m_{2,2}} + 
b_{2,2} \frac{\widetilde m_{2,1}}{\widetilde m_{2,2}} - 
b_{1,0} \frac{\widetilde m_{0,2}}{\widetilde m_{2,2}} - b_{1,1} \frac{\widetilde m_{1,2}}{\widetilde m_{2,2}},
\end{equation}
\begin{equation}
\label{f3_44__6_circle}
b_{1,1} = b_{2,0} \frac{\widetilde m_{1,1}}{2\widetilde m_{2,2}} + b_{2,1} \frac{\widetilde m_{2,1}}{\widetilde m_{2,2}} - 
b_{1,0} \frac{\widetilde m_{1,2}}{\widetilde m_{2,2}}.
\end{equation}
Substitute the expression for $b_{0,1}$ from~(\ref{f3_44__4_circle}) into~(\ref{f3_44__3_circle}) to get
\begin{equation}
\label{f3_44__3_prime}
b_{0,2} = b_{2,0} \left( \frac{\widetilde m_{0,0}}{\widetilde m_{2,2}} - \frac{\widetilde m_{1,2} \widetilde m_{1,0}}{2\widetilde m_{2,2}^2}
\right)
+ b_{2,1} \left( \frac{\widetilde m_{1,0}}{\widetilde m_{2,2}} - \frac{\widetilde m_{1,2} \widetilde m_{2,0}}{\widetilde m_{2,2}^2}
\right) +
b_{2,2} \frac{\widetilde m_{2,0}}{\widetilde m_{2,2}} + 
b_{0,0} \Delta_1,
\end{equation}
where
$$ \Delta_1 := -\frac{\widetilde m_{0,2}}{\widetilde m_{2,2}} +  
\frac{ \widetilde m_{1,2}^2 }{ 2\widetilde m_{2,2}^2 }. $$
Substitute for $b_{1,1}$ from~(\ref{f3_44__6_circle}) into~(\ref{f3_44__5_circle}) to get
\begin{equation}
\label{f3_44__5_prime}
b_{1,2} = b_{2,0} \left( \frac{\widetilde m_{0,1}}{\widetilde m_{2,2}} - \frac{\widetilde m_{1,2} \widetilde m_{1,1}}{2\widetilde m_{2,2}^2}
\right)
+ b_{2,1} \left( \frac{\widetilde m_{1,1}}{\widetilde m_{2,2}} -
\frac{\widetilde m_{1,2} \widetilde m_{2,1}}{\widetilde m_{2,2}^2}
\right)
+ 
b_{2,2} \frac{\widetilde m_{2,1}}{\widetilde m_{2,2}} + 
b_{1,0} \Delta_1.
\end{equation}
Observe that equations~(\ref{f3_44__1})-(\ref{f3_44__6}) are equivalent to equations~(\ref{f3_44__1}), (\ref{f3_44__2}),
(\ref{f3_44__3_prime}), (\ref{f3_44__4_circle}), (\ref{f3_44__5_prime}), (\ref{f3_44__6_circle}).
Notice that elements of row~$0$ are expressed in terms of $b_{0,0}$ and elements of row~$2$, while
elements of row $1$ are expressed in terms of $b_{1,0}$ and elements of row~$2$.

Substitute for $b_{0,1}$,$b_{0,2}$,$b_{1,1}$,$b_{1,2}$ from~(\ref{f3_44__3_prime}),(\ref{f3_44__4_circle}),(\ref{f3_44__5_prime}),(\ref{f3_44__6_circle}) 
into equation~(\ref{f3_44__1}). After simplifications we shall get:
$$ b_{1,0} \left( \widetilde m_{0,0} - \frac{ \widetilde m_{1,2} \widetilde m_{1,0} }{ 2\widetilde m_{2,2} } + \widetilde m_{2,0} \Delta_1
\right) =
b_{0,0} \left( \widetilde m_{0,1} - \frac{ \widetilde m_{1,2} \widetilde m_{1,1} }{ 2\widetilde m_{2,2} } + \widetilde m_{2,1} \Delta_1
\right) + $$
\begin{equation}
\label{f3_44__1_square}
+ b_{2,0} \left( \frac{\widetilde m_{2,1}\widetilde m_{0,0} - \widetilde m_{2,0}\widetilde m_{0,1}}{\widetilde m_{2,2}} +
\frac{\widetilde m_{1,2}\widetilde m_{2,0}\widetilde m_{1,1} - \widetilde m_{1,2}\widetilde m_{2,1}\widetilde m_{1,0}}{2\widetilde m_{2,2}^2}
\right).
\end{equation}
Substitute the same expressions into equation~(\ref{f3_44__2}) to get
$$ b_{1,0} \left( \widetilde m_{1,0} - \frac{ \widetilde m_{2,0} \widetilde m_{1,2} }{ \widetilde m_{2,2} }
\right) =
b_{0,0} \left( \widetilde m_{1,1} - \frac{ \widetilde m_{2,1} \widetilde m_{1,2} }{ \widetilde m_{2,2} }
\right) + $$
\begin{equation}
\label{f3_44__2_square}
+ b_{2,0} \left( \frac{\widetilde m_{2,1}\widetilde m_{1,0} - \widetilde m_{2,0}\widetilde m_{1,1}}{\widetilde m_{2,2}} 
\right).
\end{equation}
Thus, equations~(\ref{f3_44__1})-(\ref{f3_44__6}) are equivalent to equations~(\ref{f3_44__1_square}), (\ref{f3_44__2_square}),
(\ref{f3_44__3_prime}), (\ref{f3_44__4_circle}), (\ref{f3_44__5_prime}), (\ref{f3_44__6_circle}). 
Equation~(\ref{f3_44__2_square}) is equivalent to the following equation:
\begin{equation}
\label{f3_44__2_square_tilde}
b_{0,0} =  b_{1,0} \frac{\widetilde m_{2,2}}{\Delta_2} \left( 
\widetilde m_{1,0} - \frac{\widetilde m_{2,0}\widetilde m_{1,2}}{\widetilde m_{2,2}} 
\right) -
b_{2,0} \frac{\widetilde m_{2,2}}{\Delta_2}  
\frac{(\widetilde m_{2,1}\widetilde m_{1,0} - \widetilde m_{2,0}\widetilde m_{1,1})}{\widetilde m_{2,2})}{\widetilde m_{2,2}},
\end{equation}
where
$$ \Delta_2 := \widetilde m_{2,2} \widetilde m_{1,1} - \widetilde m_{1,2}^2 > 0. $$
Finally, substitute for $b_{0,0}$ from~(\ref{f3_44__2_square_tilde}) into~(\ref{f3_44__1_square}) to get
\begin{equation}
\label{f3_44__1_c}
b_{1,0} c_{1,0} =  b_{2,0} c_{2,0},
\end{equation}
where
$$ c_{1,0} := \widetilde m_{2,2} \widetilde m_{0,0} \Delta_2 - \frac{1}{2} \widetilde m_{1,2} \widetilde m_{1,0} \Delta_2 +
\widetilde m_{2,0} \widetilde m_{2,2} \Delta_1 \Delta_2 - $$
$$ - ( \widetilde m_{2,2} \widetilde m_{0,1} - \frac{1}{2} \widetilde m_{1,2} \widetilde m_{1,1} +
\widetilde m_{2,2} \widetilde m_{2,1} \Delta_1 ) ( \widetilde m_{2,2} \widetilde m_{1,0} - \widetilde m_{2,0} \widetilde m_{1,2} ), $$
$$ c_{2,0} := \Delta_2 ( \widetilde m_{2,1} \widetilde m_{0,0} - \widetilde m_{2,0} \widetilde m_{0,1} ) +
\frac{ ( \widetilde m_{1,2} \widetilde m_{2,0} \widetilde m_{1,1} - \widetilde m_{1,2} \widetilde m_{2,1} \widetilde m_{1,0} ) }
{2 \widetilde m_{2,2}} \Delta_2 - $$
$$ - ( \widetilde m_{2,1} \widetilde m_{1,0} - \widetilde m_{2,0} \widetilde m_{1,1} )
( \widetilde m_{2,2} \widetilde m_{0,1} - \frac{1}{2} \widetilde m_{1,2} \widetilde m_{1,1} 
+ \widetilde m_{2,2} \widetilde m_{2,1} \Delta_1 ). $$
After simplifications we see that there appear $\det M_0(z)$:
$$ c_{2,0} = \widetilde m_{2,1} \det M_0(z),\quad c_{1,0} = \widetilde m_{2,2} \det M_0(z). $$
Consequently, equation~(\ref{f3_44__1_c}) is equivalent to the following equation:
\begin{equation}
\label{f3_44__1_f}
b_{1,0} =  b_{2,0} \frac{ \widetilde m_{2,1} }{ \widetilde m_{2,2} }.
\end{equation}
Thus, the original system~~(\ref{f3_44__1})-(\ref{f3_44__6}) is equivalent to equations~(\ref{f3_44__1_f}), (\ref{f3_44__2_square_tilde}),
(\ref{f3_44__3_prime}), (\ref{f3_44__4_circle}), (\ref{f3_44__5_prime}), (\ref{f3_44__6_circle}). 
We can choose arbitrary real continuous on $K$ functions $b_{2,0},b_{2,1},b_{2,2}$. The rest of $b_{i,j}$ can be found from
equations~(\ref{f3_44__1_f}), (\ref{f3_44__2_square_tilde}),
(\ref{f3_44__3_prime}), (\ref{f3_44__4_circle}), (\ref{f3_44__5_prime}), (\ref{f3_44__6_circle}). 

In particular, choose $b_{2,0} = 2 \widetilde m_{2,2}^2$, $b_{2,1} = b_{2,2} = 0$, to get
$$ B(z) = $$
$$ \left(
\begin{array}{ccc}
2 \widetilde m_{2,0} \widetilde m_{2,2} & \widetilde m_{1,0}\widetilde m_{2,2} - \widetilde m_{2,0}\widetilde m_{1,2} &
2 \widetilde m_{2,2}\widetilde m_{0,0} - \widetilde m_{1,2}\widetilde m_{1,0} + \frac{\widetilde m_{2,0}\widetilde m_{1,2}^2}{\widetilde m_{2,2}}
- 2\widetilde m_{2,0}^2\\
2 \widetilde m_{2,1} \widetilde m_{2,2} & \widetilde m_{1,1}\widetilde m_{2,2} - \widetilde m_{2,1}^2 &
2 \widetilde m_{2,2}\widetilde m_{0,1} - \widetilde m_{1,2}\widetilde m_{1,1} + \frac{\widetilde m_{1,2}^3}{\widetilde m_{2,2}}
- 2\widetilde m_{2,1}\widetilde m_{0,2}\\
2\widetilde m_{2,2}^2 & 0 & 0\end{array}
\right). $$
It is easily checked that $\det B(z) \not= 0$, $z\in K$.
The proof is complete. $\Box$

\noindent
\textbf{Example 1.} Suppose that $\rho = 1$, $K = [a,b]$, $-\infty < a < b < +\infty$,
$$ M_0(z) = 
\left(
\begin{array}{cc}
1 & 0 \\
0 & \widetilde m_{1,1}(z) \end{array}                  
\right), $$
where $\widetilde m_{1,1}$ is a positive continuous function on $K$, and
$\mu$ is an arbitrary finite measure on $\mathfrak{B}(K)$. Notice that the matrix measure $M$ is now given by
$$ M(\delta) = \int_\delta M_0(z) d\mu,\qquad \delta\in\mathfrak{B}(K). $$
In this case we have (see definitions in the above proof of Theorem~\ref{t3_1}):
$$ B(z) =
\left(
\begin{array}{cc}
0 & 1 \\
\widetilde m_{1,1}(z) & 0 \end{array}                  
\right), $$
$$ B^{-1}(z) =
\left(
\begin{array}{cc}
0 & \frac{1}{\widetilde m_{1,1}(z)} \\
1 & 0 \end{array}                  
\right),\quad
C(z) =
\left(
\begin{array}{cc}
0 & z \\
z \widetilde m_{1,1}(z) & \widetilde m_{1,1}(z) \end{array}                  
\right). $$
Notice that
$$ B^{-1}(z) C(z) = A(z) =
\left(
\begin{array}{cc}
z & 1 \\
0 & z \end{array}                  
\right), $$
but
$$ C(z) B^{-1}(z) =
\left(
\begin{array}{cc}
z & 0 \\
\widetilde m_{1,1}(z) & z \end{array}                  
\right). $$
We have
$$ B(z) M_0(z) =
\left(
\begin{array}{cc}
0 & \widetilde m_{1,1}(z) \\
\widetilde m_{1,1}(z) & 0 \end{array}                  
\right),\quad  
C(z) M_0(z) =
\left(
\begin{array}{cc}
0 & z \widetilde m_{1,1}(z) \\
z \widetilde m_{1,1}(z) & \widetilde m_{1,1}^2(z) \end{array}                  
\right), $$
and
$$ B^{-1}(z) M_0(z) =
\left(
\begin{array}{cc}
0 & 1 \\
1 & 0 \end{array}                       
\right). $$

\begin{corollary}
\label{c3_1}
In conditions of Theorem~\ref{t3_1}, suppose additionally that condition~(\ref{f2_7}) holds 
for arbitrary non-zero $p\in\mathbb{P}$.
Let
\begin{equation}
\label{f3_50}
\sigma(f,g) = [f,g]_{\mathcal{B}^{-1}} :=
(\mathcal{B}^{-1} f,g)_{L^2(M)},\qquad f,g\in C(K).
\end{equation}
Then the associated Sobolev orthogonal polynomials $\{ y_n(z) \}_{n=0}^\infty$, as in~(\ref{f2_10}), satisfy the following
relations:
\begin{equation}
\label{f3_54}
\left[
\mathcal{B} y_n(\mathcal{A}) 1, y_m(\mathcal{A}) 1 
\right]_{\mathcal{B}^{-1}} = A_n \delta_{n,m},\qquad n,m\in\mathbb{Z}_+.
\end{equation}
Operators $\mathcal{A},\mathcal{B}$ are symmetric with respect to $\sigma$:
$$ \sigma(\mathcal{A}f,g) = \sigma(f,\mathcal{A}g),\quad \sigma(\mathcal{B}u,v) 
= \sigma(u,\mathcal{B}v), $$
for all $f,g\in A^2(M)$, $u,v\in C(K)$
\end{corollary}

\noindent
\textbf{Proof.} The proof is straightforward. $\Box$

It looks natural to state the following conjecture.

\noindent
\textbf{Conjecture 1.}
Theorem~\ref{t3_1} is true, if we replace the assumption $\rho\in\{ 1,2 \}$ by 
$\rho\in\mathbb{N}$. 

As we have seen in the proof of Theorem~\ref{t3_1}, there already appeared huge expressions for $\rho=2$.
On the other hand, after proper sorting of the equations and lots of simplifications there appeared $\det M_0(z)$. It looks promising for
the validity of the general case.

}
\end{document}